\documentclass[11pt, leqno]{article} 
\usepackage{a4}
\usepackage{amssymb, amsmath, amsfonts,amsthm, epsfig, graphics, tikz, float}

\begin{document}
\title{On a preconditioner for time domain boundary element methods}
\author{ Heiko Gimperlein\thanks{Maxwell Institute for Mathematical Sciences and Department of Mathematics, Heriot--Watt University, Edinburgh, EH14 4AS, United Kingdom, email: h.gimperlein@hw.ac.uk, ds221@hw.ac.uk}  \thanks{Institute for Mathematics, University of Paderborn, Warburger Str.~100, 33098 Paderborn, Germany. \newline We thank Dugald Duncan for discussions around extrapolation as a possible preconditioner. H.~G.~acknowledges support by ERC Advanced Grant HARG 268105 and the EPSRC Impact Acceleration Account. } \and David Stark${}^\ast$ }
\date{}

\maketitle \vskip 0.5cm
\begin{abstract}
\noindent We  propose a time stepping scheme for the space-time systems obtained from Galerkin time-domain boundary element methods for the wave equation. Based on extrapolation, the method proves stable, becomes exact for increasing degrees of freedom and can be used either as a preconditioner, or as an efficient standalone solver for scattering problems with smooth solutions. It also significantly reduces the number of GMRES iterations for screen problems, with less regularity, and we explore its limitations for enriched methods based on non-polynomial approximation spaces. 
\end{abstract}

\section{Introduction}

Efficient and accurate computational methods to simulate sound emission in space and time are of interest from the modeling of environmental noise to acoustic scattering \cite{ajrt, Banz, sayas, mic}. Computations in time domain are of particular interest for problems beyond the reach of frequency domain methods, such as the simulation of transient dynamics, moving sound sources or nonlinear and dynamical contact problems. 

Galerkin time domain boundary element methods prove to be stable and accurate in long--time computations and are competitive with frequency domain methods for realistic problems \cite{survey}. Practical implementations, however, are based on marching-on-in-time (MOT) time stepping schemes, corresponding to a Petrov-Galerkin method with piecewise constant test functions or Dirac distributions in time.  Relevant works on the numerical implementation include the Ph.D.~thesis of Terrasse \cite{terrasse93} and \cite{dfhd}, which made the methods competitive for commercial applications. 

In practice, the different choices of ansatz and test functions for MOT schemes may lead to instabilities. Proper Galerkin methods are not only provably stable, but they have also attracted interest from at least three different perspectives: Rigorous a posteriori error estimates give rise to efficient adaptive mesh refinement procedures \cite{Gimperlein, Glafke, GM1, GM2, sv, veit}; non-polynomial basis functions and efficient assembly of the algebraic system \cite{sv2, veit, czech}; formulations based on the physical energy \cite{aimineu, aiminumeralg, aimi5}. In this context, efficient preconditioners have been of current  interest. These will also be crucial for efficient higher order methods.

This work studies a time stepping scheme to solve the space time systems arising from time domain Galerkin boundary element discretisations for higher order test functions. It approximates the algebraic system using extrapolation and may be used as either a preconditioner for the full space-time equation or as a fast, quasi-exact stand-alone solver of the integral equation for large numbers of degrees of freedom. Our method inherits the approximation properties and long-time stability of the Galerkin method and has been used, but not thoroughly documented in recent works \cite{Gimperlein3, Gimperlein, pupaper}. The rigorous numerical analysis of the surprisingly good stability properties of the proposed method for the standard $h$-TDBEM under mesh refinements remains open.\\

\noindent \emph{Outline of the article:} In Section 2 we recall the boundary integral formulation of the wave equation and its numerical discretization using boundary elements. The MOT time stepping scheme for the resulting space-time system is described in Section 3. Section 4 investigates the scheme as a preconditioner and as a standalone solver in numerical experiments on closed surfaces, screens and with non-polynomial basis functions.

\section{Problem formulation}

We consider transient sound radiation problems in the exterior of a scatterer $\Omega^{-}$, where $\Omega^{-}$ is a bounded polygon or a screen with connected complement 
$\Omega = \mathbb{R}^3\setminus \overline{\Omega^{-}}$. The acoustic sound pressure field $u$ due to an incident field or sources on $\Gamma = \partial \Omega$ satisfies the linear wave equation for $t \in \mathbb{R}$:
 \begin{align}\label{waveeq}
c^{-2}\partial_t^2 u(t,{\bf x}) - \Delta u(t,{\bf x}) &= 0 \quad \text{for $(t, \mathbf{x})\in\mathbb{R}\times\Omega$},\\ \nonumber u(t,{\bf x})=f(t,{\bf x}) \quad \text{for ${\bf x}\in\Gamma$}, &\quad u(t,{\bf x}) = 0 \quad \text{for $t \leq 0$.}
 \end{align}
Here $c$ is the wave speed, and in the following we set $c=1$ for  simplicity. A single-layer ansatz for $u$,
\begin{equation}\label{singlay}
u(t,{\bf x}) = {\int_\Gamma \frac{\phi (t-|{\bf x}-{\bf y}|,{\bf y})}{2\pi|{\bf x}-{\bf y}|}ds_y},
\end{equation}
results in an equivalent weak formulation of \eqref{waveeq} as an integral equation of the first kind in space-time anisotropic Sobolev spaces \cite{survey, Nezhi}:\\ {Find $\phi \in H^1_{\sigma}(\mathbb{R}^+,\widetilde{H}^{-\frac{1}{2}}(\Gamma))$ such that} for all $\psi\in H^1_{\sigma}(\mathbb{R}^+,\widetilde{H}^{-\frac{1}{2}}(\Gamma))$
\begin{equation}\label{weakform}
 \int_0^\infty\int_\Gamma (V \phi(t,{\bf x})) \partial_t\psi(t,{\bf x}) \ ds_x\ d_\sigma t = \int_0^\infty \int_\Gamma f(t,{\bf x}) \partial_t\psi(t,{\bf x}) \ ds_x\ d_\sigma t\ , 
\end{equation}
where 
$$V \phi(t,{\bf x}) = \int_\Gamma \frac{\phi(t-|{\bf x}-{\bf y}|,{\bf y})}{2\pi|{\bf x}-{\bf y}|} ds_y\ ,$$
and $d_\sigma t = e^{-2 \sigma t}dt$. A theoretical analysis requires $\sigma>0$, but practical computations use $\sigma=0$ \cite{bamberger,survey}.\\

We study time dependent boundary element methods to solve \eqref{weakform}, based on approximations {by} piecewise polynomial ansatz and test functions from the space $V^{p,q}_{h,\Delta t}$ spanned by
\begin{equation}\label{testfunc}
{\phi_{mi}(t,{\bf x}) = \widetilde{\Lambda}_m(t) {\Lambda}_i({\bf x})}\ .
\end{equation}
Here, ${\Lambda}_i$ a piecewise polynomial shape function of degree $p$ in space and $\widetilde{\Lambda}_m$ a corresponding shape function of degree $q$ in time. For $p\geq 1$, resp.~$q\geq 1$, the shape functions are assumed to be continuous. 


We obtain a numerical scheme for the weak formulation \eqref{weakform}: \\Find $\phi_{h,\Delta t}\in V^{p,q}_{h,\Delta t}$ such that for all $\psi_{h,\Delta t} \in V^{p,q}_{h,\Delta t}$
\begin{equation}\label{weakhform}
 {\int_0^\infty\int_\Gamma \big(V \phi_{h,\Delta t}(t,{\bf x})\big) \partial_t\psi_{h,\Delta t}(t,{\bf x}) \ ds_x\ dt = \int_0^\infty \int_\Gamma f(t,{\bf x}) \partial_t\psi_{h,\Delta t}(t,{\bf x}) \ ds_x\ dt}\ .
\end{equation}
From $\phi_{h,\Delta t}$, the sound pressure $u_{h,\Delta t}$ is obtained in $\Omega$ by evaluating the integral in \eqref{singlay} numerically. \\

\begin{figure}[htbp]
 \centering
 \includegraphics[height=3.2cm]{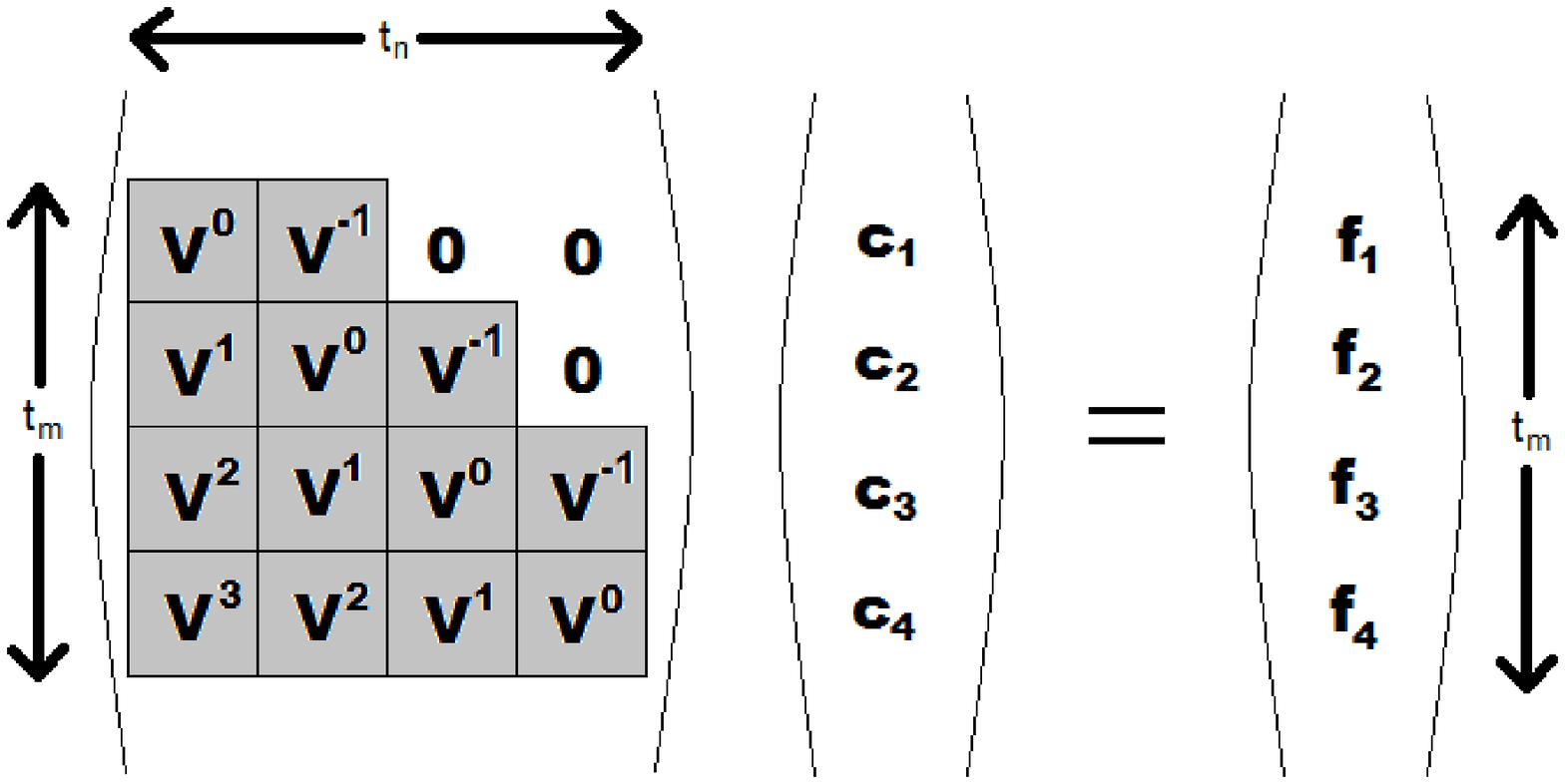}
 \caption{Full space-time system showing the bands $\mathbf{V}^j$.}
 \label{system}
\end{figure}

{With $\phi_{h,\Delta t}(t,{\bf x}) = \sum_{m}\sum_i {\bf c}_m^i \phi_{mi}(t,{\bf x})$, equation \eqref{weakhform} leads to a linear system ${\bf V} {\bf c} = {\bf F}$ of equations for the coefficients ${\bf c}_m^i$ in space-time, Figure \ref{system}.} 
{Here the stiffness matrix ${\bf V}$ has a block--lower--Hessenberg form with blocks
$${\bf V}^{m-n}_{ij} = \int_0^\infty \int_\Gamma \big(V \widetilde{\Lambda}_m(t){\Lambda}_i({\bf x})\big) \partial_t\widetilde{\Lambda}_n(t){\Lambda}_j({\bf x}) ds_x dt\ ,$$
corresponding to the time steps. The vector on the right hand side  is given by $${\bf f}_n^j =\int_0^\infty\int_\Gamma f(t,{\bf x}) \partial_t\widetilde{\Lambda}_n(t){\Lambda}_j({\bf x}) \ ds_x\ dt \ .$$}

A main challenge is the accurate assembly of $\mathbf{V}$. After
an analytical evaluation of the time integral, the $y$ integral requires integration over geometrically complicated intersections of triangles with light cone shells, with a singular integrand $|{\bf x}-{\bf y}|^{-1}$.  It is evaluated in polar coordinates with a geometrically-graded $hp$-composite Gauss quadrature \cite{survey}. A regular Gauss quadrature is used for the $x$ integral.\\

\begin{figure}[htbp]
 \centering
 \includegraphics[height=3.2cm]{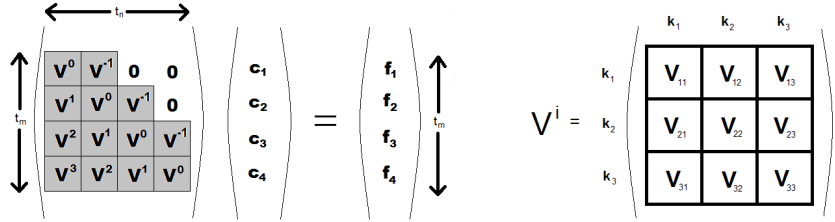}
 \caption{Full space-time system for PU TDBEM.}
 \label{system2}
\end{figure}

We also consider a time domain \emph{partition of unity} method.  Here, the ansatz and test functions are, instead of \eqref{testfunc}, given by travelling plane waves {\begin{equation}{\phi_{mi}}(t,{\bf x})=\widetilde{\Lambda}_m(t)\Lambda_i({\bf x})\cos(\omega_{i} t - {\bf k}_i\cdot{\bf x} + \sigma_i )\ ,\label{pufunc}\end{equation} where $\omega_i=\|{\bf k}_i\|_2$, and $\sigma_i=\{0,\frac{\pi}{2}\}$. 
The numerical scheme for the weak formulation \eqref{weakform} reads as follows: \\Find $\phi_{h,\Delta t}(t,{\bf x}) = \sum_{m}\sum_i {\bf c}_m^i \phi_{mi}(t,{\bf x})$ such that 
\begin{equation}\label{weakhformpu}
 {\int_0^\infty\int_\Gamma \big(V \phi_{h,\Delta t}(t,{\bf x})\big) \partial_t\psi_{h,\Delta t}(t,{\bf x}) \ ds_x\ dt = \int_0^\infty \int_\Gamma f(t,{\bf x}) \partial_t\psi_{h,\Delta t}(t,{\bf x}) \ ds_x\ dt}
\end{equation}
for all $\psi_{h,\Delta t}(t,{\bf x}) = \sum_{n}\sum_j {\bf c}_n^j \phi_{nj}(t,{\bf x})$.
This again leads to a linear system ${\bf V} {\bf c} = {\bf F}$ of equations for the coefficients ${\bf c}_m^i$ in space-time, Figure \ref{system2}, and in this case each of the time step blocks ${\bf V}^{m-n}$ in addition decomposes into blocks for the individual ${\bf k}_i$. The blocks of the stiffness matrix ${\bf V}$ are given by 
$${\bf V}^{m-n}_{ij} = \int_0^\infty \int_\Gamma \big(V \phi_{mi}(t,{\bf x})\big) \partial_t\phi_{nj}(t,{\bf x}) ds_x dt\ ,$$
corresponding to the time steps, and the vector on the right hand side  is $${\bf f}_n^j =\int_0^\infty\int_\Gamma f(t,{\bf x}) \partial_t\phi_{nj}(t,{\bf x}) \ ds_x\ dt \ .$$}

\section{Time stepping scheme for the space-time system}

Because the blocks $\mathbf{V}^{-1}$ above the diagonal are nonzero for test functions of degree $\geq 1$ in time, in general the full space-time system \eqref{weakhform} has to be solved for all time steps at once. Solving such large systems iteratively, e.g.~using GMRES, can take a substantial number of iterations (and time) to achieve an acceptable tolerance. On the other hand, for test functions of degree $0$ in time, $\mathbf{V}^{-1}=\textbf{0}$, and backsubstitution leads to an efficient and popular time stepping scheme (MOT, \emph{marching in on time}). 

In this note we compare \eqref{weakhform} with a similar time stepping scheme obtained from a lower-triangular approximation $\widetilde{\mathbf{V}}$ for $\mathbf{V}$. We use $\widetilde{\mathbf{V}}$ both as a preconditioner for GMRES and investigate its accuracy as a stand-alone solver. \\

For simplicity we restrict ourselves to piecewise linear functions in time, $q=1$. The matrix $\widetilde{\mathbf{V}}$ is obtained from a linear extrapolation in time to approximate the solution vector for the next time step: $$\mathbf{c}_{j+1} \simeq \mathbf{c}_{j}+(\mathbf{c}_{j}-\mathbf{c}_{j-1})\ .$$
Inserting this into the $j$-th equation in Figure 1,
\begin{align*}
{\mathbf{f}_j}&=\mathbf{V}^{j-1} \mathbf{c}_{1}+\cdots +\mathbf{V}^2 \mathbf{c}_{j-2}+\mathbf{V}^1 \mathbf{c}_{j-1}+\mathbf{V}^0\mathbf{c}_j+\mathbf{V}^{-1}\mathbf{c}_{j+1}\ ,
\end{align*}
we obtain
\begin{align*}
{\mathbf{f}_j} &\simeq \mathbf{V}^{j-1} \mathbf{c}_{1}+\cdots +\mathbf{V}^2 \mathbf{c}_{j-2}+\mathbf{V}^1 \mathbf{c}_{j-1}+\mathbf{V}^0\mathbf{c}_j+\mathbf{V}^{-1}(2\mathbf{c}_j-\mathbf{c}_{j-1})\\
& =	\mathbf{V}^{j-1} \mathbf{c}_{1}+\cdots +\mathbf{V}^2 \mathbf{c}_{j-2}+\underbrace{(\mathbf{V}^1 -\mathbf{V}^{-1})}_{\widetilde{{\bf V}}^1}\mathbf{c}_{j-1}+\underbrace{(\mathbf{V}^0 +2\mathbf{V}^{-1})}_{\widetilde{{\bf V}}^0}\mathbf{c}_j \ .
\end{align*}
Therefore, we define $\widetilde{\mathbf{V}}$ as the lower-triangular block Toeplitz matrix with bands $\widetilde{\mathbf{V}}^i = {\mathbf{V}}^i$ for $i>1$, $\widetilde{\mathbf{V}}^1=\mathbf{V}^1-\mathbf{V}^{-1}$ and $\widetilde{\mathbf{V}}^0=\mathbf{V}^0 +2\mathbf{V}^{-1}$. Solving the lower triangular system $\widetilde{\mathbf{V}}\mathbf{c}=\mathbf{F}$ by backsubstitution, we obtain a modified MOT scheme.

\section{Numerical results}

We explore the relevance of $\widetilde{\mathbf{V}}$ {for the $h$-version TDBEM} in three model problems for plane-wave scattering from an icosahedron or sphere, resp.~a screen. For all three problems, we use ansatz and test functions in $V^{0,1}_{h,\Delta t}$ {and solve the block-Hessenberg system \eqref{weakhform}. A fourth example discusses the relevance of $\widetilde{\mathbf{V}}$ for the partition of unity TDBEM, with block-Hessenberg system \eqref{weakhformpu}.} 

As a preconditioner within GMRES, the computational cost is similar to a matrix-vector multiplication and approximately doubles the cost of each GMRES iteration.\\

\noindent \textbf{Example 1:} First we consider a sphere $\Gamma = S^2$ with right hand side $f(t,x) = \sin(2t)^5$. The exact solution is given by $\phi(t,x) = 10 \cos(2t) \sin(2t)^4$.  We choose a time interval $[0,T]=[0,2.5]$ and uniform meshes of 320 triangles with time step $\Delta t=0.04$, resp.~$1280$ triangles with $\Delta t=0.02$, thereby keeping the CFL number $\frac{\Delta t}{\Delta x}$ fixed.

\begin{figure}[H]
  \centering
   \includegraphics[width=12cm]{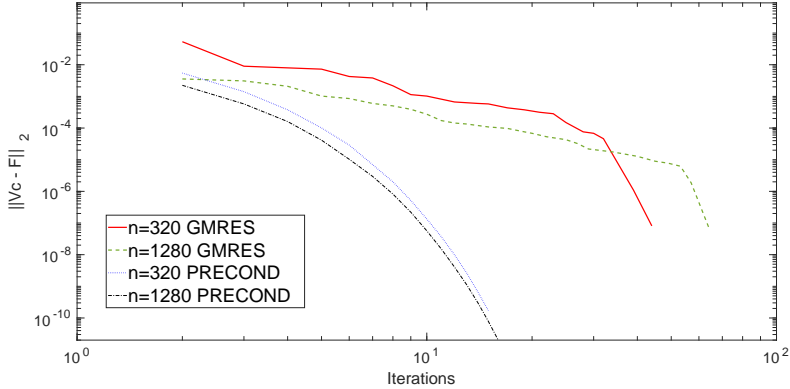}
   \caption{Residual of GMRES for $\Gamma = S^2$ with $320$, $1280$ triangles.}\label{sphereresid}
\end{figure}

We solve the resulting system with standard GMRES and with the  $\widetilde{\mathbf{V}}$-preconditioned GMRES until the standard residual error indicator for GMRES is smaller than $10^{-9}$; for GMRES this indicator is $\|{\bf V\bf c-\bf F}\|_2$, while it is $\|\widetilde{\bf V}^{-1}({\bf V\bf c-\bf F})\|_2$  for the preconditioned GMRES. {Figure \ref{sphereresid}} compares the residual $\|{\bf V\bf c-\bf F}\|_2$ in each iteration. Note that the preconditioner not only reduces the number of iterations compared to the standard GMRES, but fewer preconditioned GMRES iterations are required for $1280$ triangles: The extrapolation in time used to construct $\widetilde{\mathbf{V}}$ becomes exact as $\Delta t$ tends to $0$, provided $\phi$ is sufficiently regular. We explore the stability and accuracy of the preconditioner in detail for more realistic scattering problems below. \\

\begin{figure}[H]
  \centering
   \includegraphics[width=12cm]{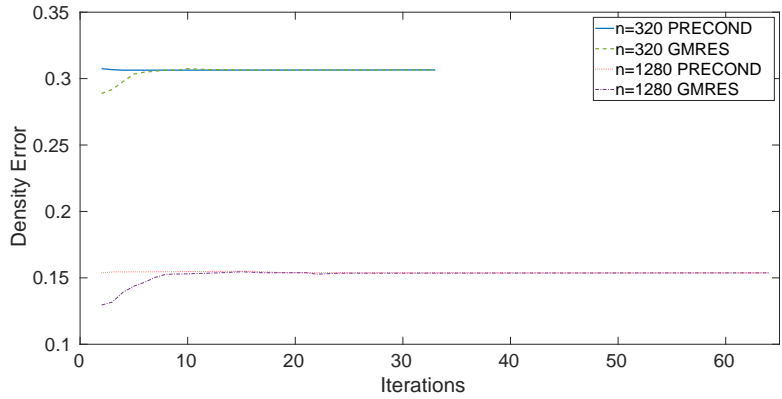}
   \caption{Density error $\|\phi - \phi_{h, \Delta t} \|_{L^2([0,T]\times \Gamma)}$ for $\Gamma = S^2$ with $320$, $1280$ triangles. }
\label{spheredens}
\end{figure}

For the current example, we can easily compare the resulting error of GMRES with the discretization error in the density $\phi$. Figure \ref{spheredens} shows that the density error $\|\phi - \phi_{h, \Delta t} \|_{L^2([0,T]\times \Gamma)}$ is dominated by the discretization error for both $320$ and $1280$ triangles. However, the contribution of the GMRES error increases for finer meshes and is nonnegligible in the first iterations. 

\begin{figure}[H]
  \centering
\label{sphereresit}
   \includegraphics[width=12cm]{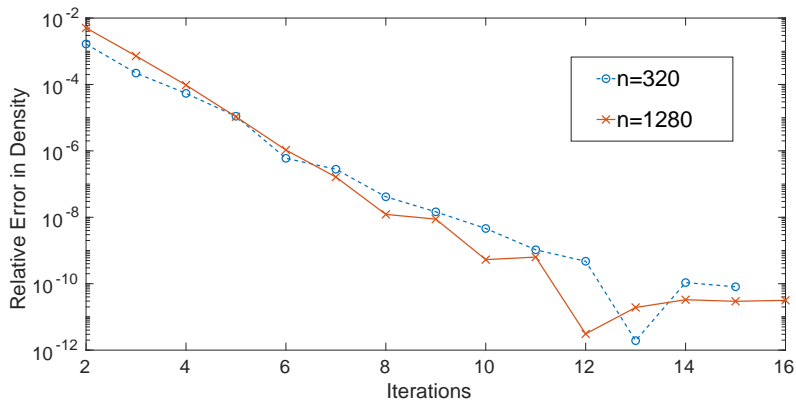}
   \caption{Relative error in Figure \ref{spheredens} for the preconditioned GMRES with 320, 1280 triangles. }\label{relerrsphere}
\end{figure}

For the preconditioned GMRES the error is negligible already after one iteration. Indeed, Figure \ref{relerrsphere} shows the relative error from the results of Figure \ref{spheredens} for the preconditioned GMRES. For both $320$ and $1280$ DOF, the density converges to the exact density of the discretized problem at approximately the same rate and magnitude.

\begin{figure}[H]
  \centering
\label{sphereresit}
   \includegraphics[width=12cm]{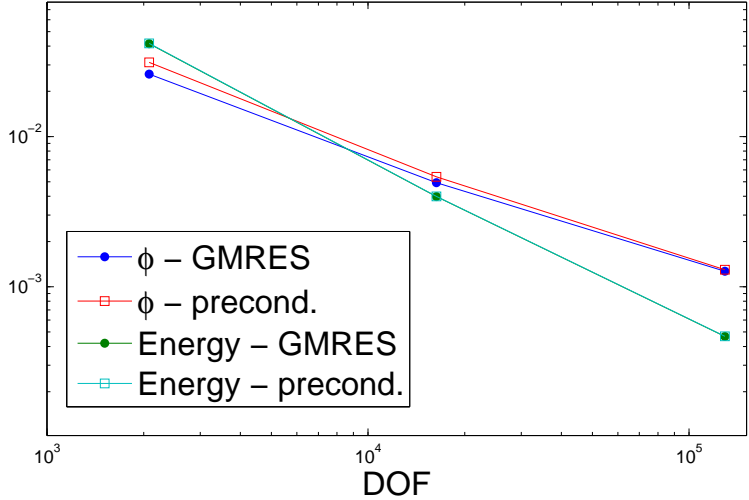}
   \caption{Convergence of TDBEM using preconditioned GMRES, resp.~the preconditioner as standalone solver: relative errors in energy and ${L^2([0,T]\times \Gamma)}$. }\label{convplot}
\end{figure}

We also consider the error between the solutions to $\widetilde{\bf V}\bf c=\bf F$  and ${\bf V}\bf c=\bf F$, i.e.~the use of the preconditioner as a stand-alone solver: As described in Table \ref{tab1step}, the relative error in {the energy $E=\frac{1}{2} {\bf c}^T {\bf V}{\bf c}-{\bf F}^T \bf c$} is well below  $1\%$ after an application of the preconditioner, resp.~a single step of preconditioned GMRES. For $1280$ triangles, a single iteration of the preconditioner yields a practically exact energy. Analogous results are obtained for the density instead of the energy, see e.g.~Figure \ref{relerrsphere}. Figure \ref{convplot} finally shows that standard convergence plots for the error of either the energy or the norm $L^2([0,T]\times \Gamma)$ may be calculated from the preconditioner as a standalone solver: The convergence rates coincide with those calculated from numerical solutions using the preconditioned GMRES solver.

\begin{table}
\begin{center}
\begin{tabular}{clccccc}
Example & Geometry &spatial DOF
    & Energy 
& Preconditioner & 1 precond. GMRES \\ \hline
1&sphere&320&8.5692&8.5470 (.26\%)&8.5481 (.25\%) \\
&&1280&8.6059&8.6059 ($\ll$ 1\%)&8.5954 (.12\%)\\
2&icosahedron&320&20.5388&21.4801 (4.6\%)&21.6785 (5.5\%)\\
&&1280&19.8796&20.1434 (1.3\%)&20.159 (1.4\%)\\
3&screen&288&0.4233&0.4497 (6.2\%)&0.4522 (6.8\%)\\
&&1250&0.4589&0.4716 (2.8\%)&0.4721 (2.9\%)\\
4&icosahedron&7 enrichments&23.6226&23.126 (2.1\%)&16.897 (28.5\%) \\
&(using PU)&15 enrichments&22.9151&22.685 (1.0\%)&20.947 (8.6\%)
\end{tabular}
\end{center}
\caption{\label{tab1step}Energy and relative errors in energy for preconditioner, resp.~a single step of preconditioned GMRES.}
\end{table}

\begin{figure}[H]
  \centering
   \includegraphics[width=10cm]{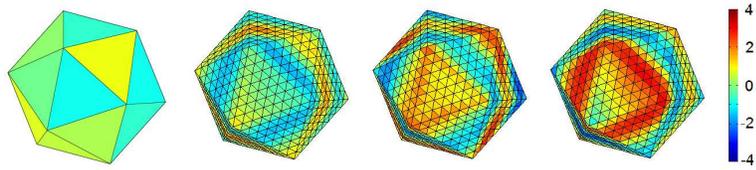}
   \caption{Meshes for icosahedron with 20 and 1280 triangles. } \label{icos}
\end{figure}

\noindent \textbf{Example 2:} The second experiment considers plane-wave scattering from an icosahedron of diameter $2$ centred in $\bf 0$. The right hand side is given by $f(t,{\bf x})=\exp(-25/t^2) \cos(\omega_f t-{\bf k}_f{\bf x})$, where ${\bf k}_f=(1.5, 3, 8.5)$ and $\omega_f= \|{\bf k}_f\|_2$. For the discretization we consider uniform meshes with $20$, $80$, $320$, $1280$, $5120$ triangles and corresponding time steps $\Delta t = 0.16,\ 0.08,\ 0.04,\ 0.02,\ 0.01$ in the time interval [0,T]=[0,5]. The icosahedron and snapshots of the numerical solution for $1280$ triangles are shown in Figure \ref{icos}.

\begin{figure}[H]
  \centering
   \includegraphics[width=12cm]{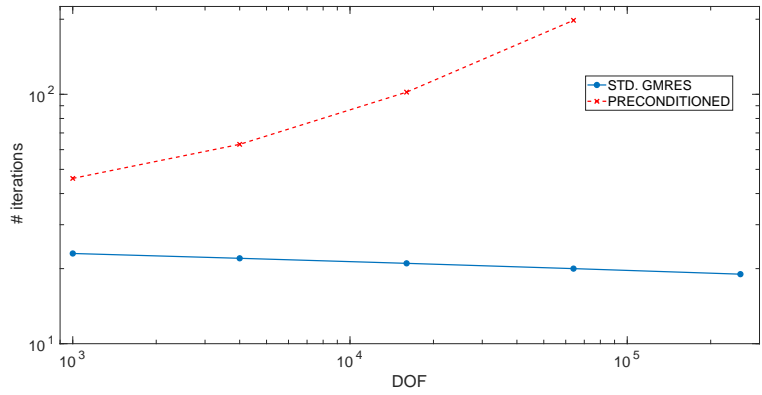}
   \caption{Iterations vs.~DOF to achieve a residual $<10^{-9}$ on the icosahedron. }\label{precondico}
\end{figure}

\begin{figure}[H]
  \centering
   \includegraphics[width=12cm]{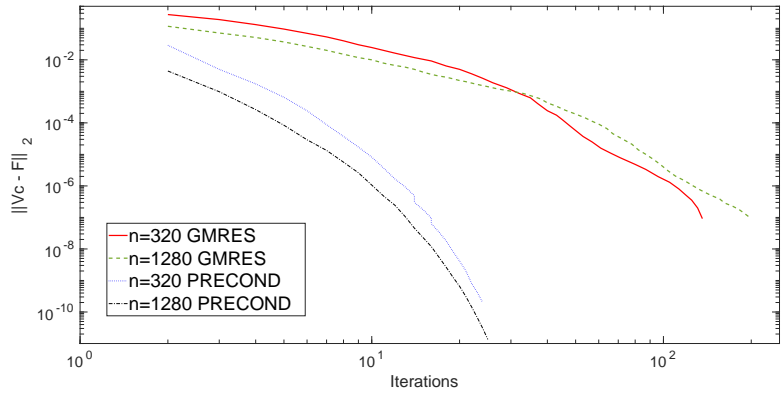}
   \caption{Residual of GMRES for icosahedron with $320$, $1280$ triangles.}\label{precondres}
\end{figure}

\begin{figure}[H]
  \centering
   \includegraphics[width=12cm]{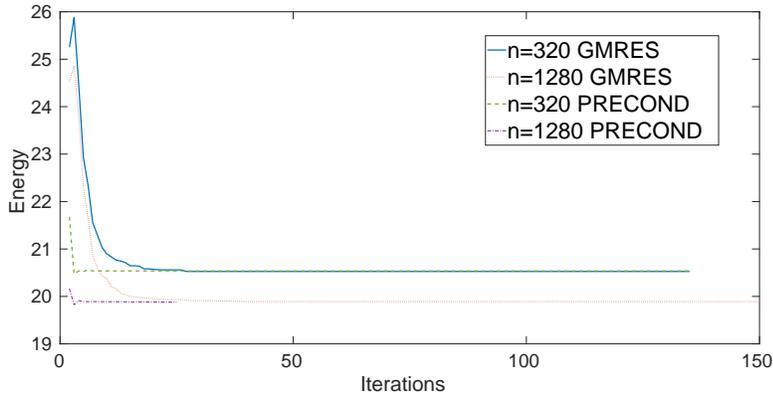}
   \caption{Energy convergence for icosahedron with $320$, $1280$ triangles.}\label{icosenergy}
\end{figure}

Figure \ref{precondico} depicts the number of iterations needed to achieve a residual error indicator as in Example 1 of $<10^{-7}$ for increasing degrees of freedom. Unlike the standard GMRES solver, the preconditioned GMRES is stable, and its iteration count even decreases from 25 to 19 for the considered meshes. Figure \ref{precondres} shows the convergence of the residuals in each iteration for $320$ and $1280$ triangles, with and without preconditioner. The preconditioner leads to a consistently larger reduction of the residual in every step, and for the finer mesh the norm of the residual of the preconditioned GMRES is reduced by a multiplicative constant. Similarly, Figure \ref{icosenergy} shows the improved convergence of the energy.
This reinforces our conclusion from Example 1 that the preconditioner becomes exact as $\Delta t$ to $0$ in a realistic scattering problem. 

Unlike in Example 1, however, Table \ref{tab1step} shows that the error of a single step of preconditioned GMRES, resp.~the preconditioner alone, does not decrease the error of the energy below $1\%$ for this more realistic scattering problem. The relative error of the energy decreases to around $1\%$, though, for $1280$ triangles, an accuracy that may suffice for large scale engineering applications.

\noindent \textbf{Example 3:} We now investigate plane-wave scattering from a screen, the square $\Gamma = [0, 0.5]^2 \times \{z=0\} \subset \mathbb{R}^3$. Here one expects edge and corner singularities and a corresponding nonsmooth behaviour of the solution in time. The resulting effect on the preconditioner is of particular interest, as it is based on a Taylor expansion in time. 

\begin{figure}[H]
  \centering
   \includegraphics[width=8cm]{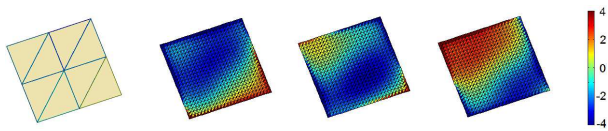}
   \caption{Meshes for screen with 8 and 1250 triangles. } \label{squaresol}
\end{figure}

We again choose a plane wave right hand side $f(t,{\bf x})=\exp(-4/t^2) \cos(\omega_f t-{\bf k}_f{\bf x})$, with ${\bf k}_f=(2, 2, 2)$ and $\omega_f = \|{\bf k}_f\|_2$. The square $\Gamma$ is discretized with uniform meshes of $288$ and $1250$ triangles, and we use corresponding time steps of $\Delta t = 0.1$, resp.~$0.05$, for the time interval $[0,T]=[0,2.5]$.  Figure \ref{squaresol} shows the mesh with $8$ triangles as well as snapshots of the numerical solution for $1250$ triangles. 

\begin{figure}[H]
  \centering
   \includegraphics[width=12cm]{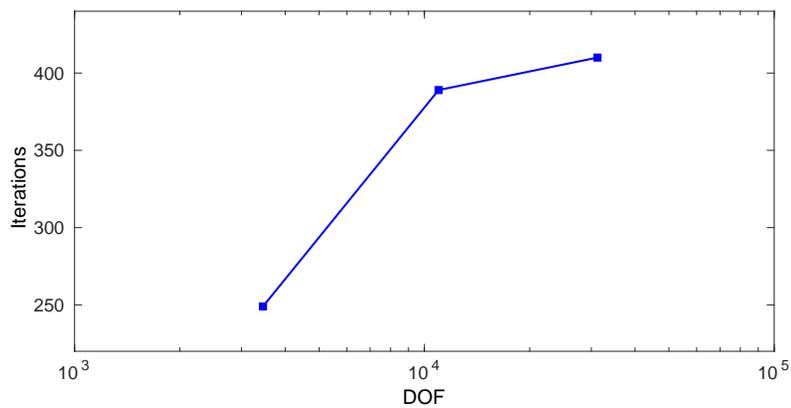}
   \caption{Iterations vs.~DOF to achieve a residual $<10^{-7}$ with preconditioned GMRES on the screen. } \label{prec_screen}
\end{figure}

In Figure \ref{prec_screen} we plot the number of iterations needed to achieve an error indicator as in Example 1 of $<10^{-7}$ for increasing degrees of freedom.  The number of iterations for the preconditioned GMRES is no longer stable, but grows with the DOF, though with a smaller exponent than without preconditioner.\\

\begin{figure}[H]
  \centering \label{precondres2}
   \includegraphics[width=12cm]{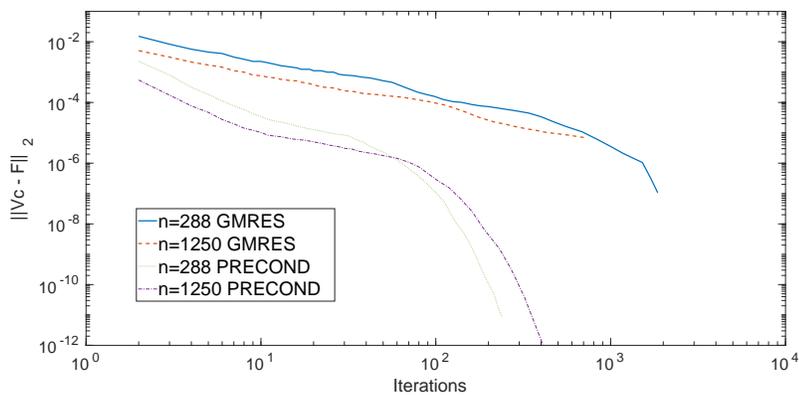}
   \caption{Residual of GMRES for screen with $288$, $1250$ triangles. }\label{resitscreen}
\end{figure}

\begin{figure}[H]
  \centering
   \includegraphics[width=12cm]{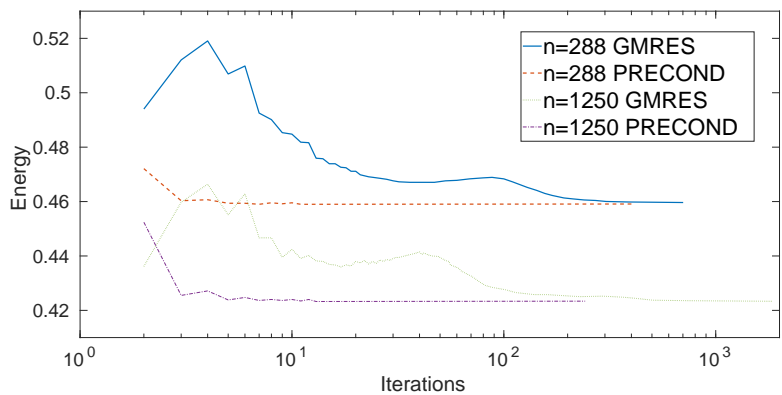}
   \caption{Energy convergence for screen with $288$, $1250$ triangles. }\label{energyscreen}
\end{figure}

 Figure \ref{resitscreen} shows the convergence of the residuals in each iteration for $288$ and $1250$ triangles, with and without preconditioner, while Figure \ref{energyscreen} depicts the convergence of the energy. The results reconfirm the relevance of the preconditioner also for the less regular solutions that naturally arise on a screen. However, unlike for the closed surfaces of Examples 1 and 2, an application of the preconditioner or a single preconditioned GMRES iteration still yields a larger error in the energy of around $6\%$ for $288$ triangles and less than $3\%$ for $1250$ triangles, see Table \ref{tab1step}.

\noindent \textbf{Example 4:} We finally solve the space-time system for a partition of unity method based on ansatz and test functions of the form \eqref{pufunc}, with piecewise linear shape functions in time and piecewise constant shape functions in space. As in Example 2 we consider an icosahedron of diameter $2$ centred in $\bf 0$. The right hand side is given by $f(t,{\bf x})=\exp(-25/t^2) \cos(\omega_f t-{\bf k}_f{\bf x})$, where ${\bf k}_f=(8.5,3,0.5)$ and $\omega_f= \|{\bf k}_f\|_2$. For the discretization we consider the $20$ faces of the icosahedron, depicted in Figure \ref{icos}, and use up to $15$ travelling plane waves per triangle. We fix the time step $\Delta t = 0.1$ in the time interval $[0,T]=[0,2.5]$. 

\begin{figure}[H]
  \centering
   \includegraphics[width=12cm]{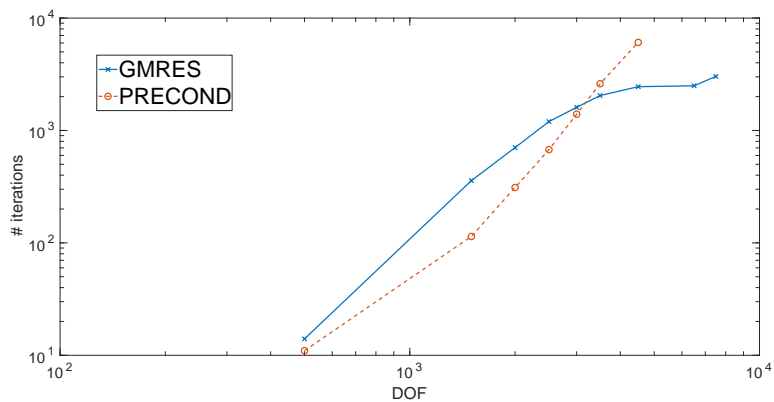}
   \caption{Iterations vs.~DOF to achieve a residual $<10^{-6}$ with a partition of unity method.}\label{precondpu}
\end{figure}

\begin{figure}[H]
  \centering
   \includegraphics[width=12cm]{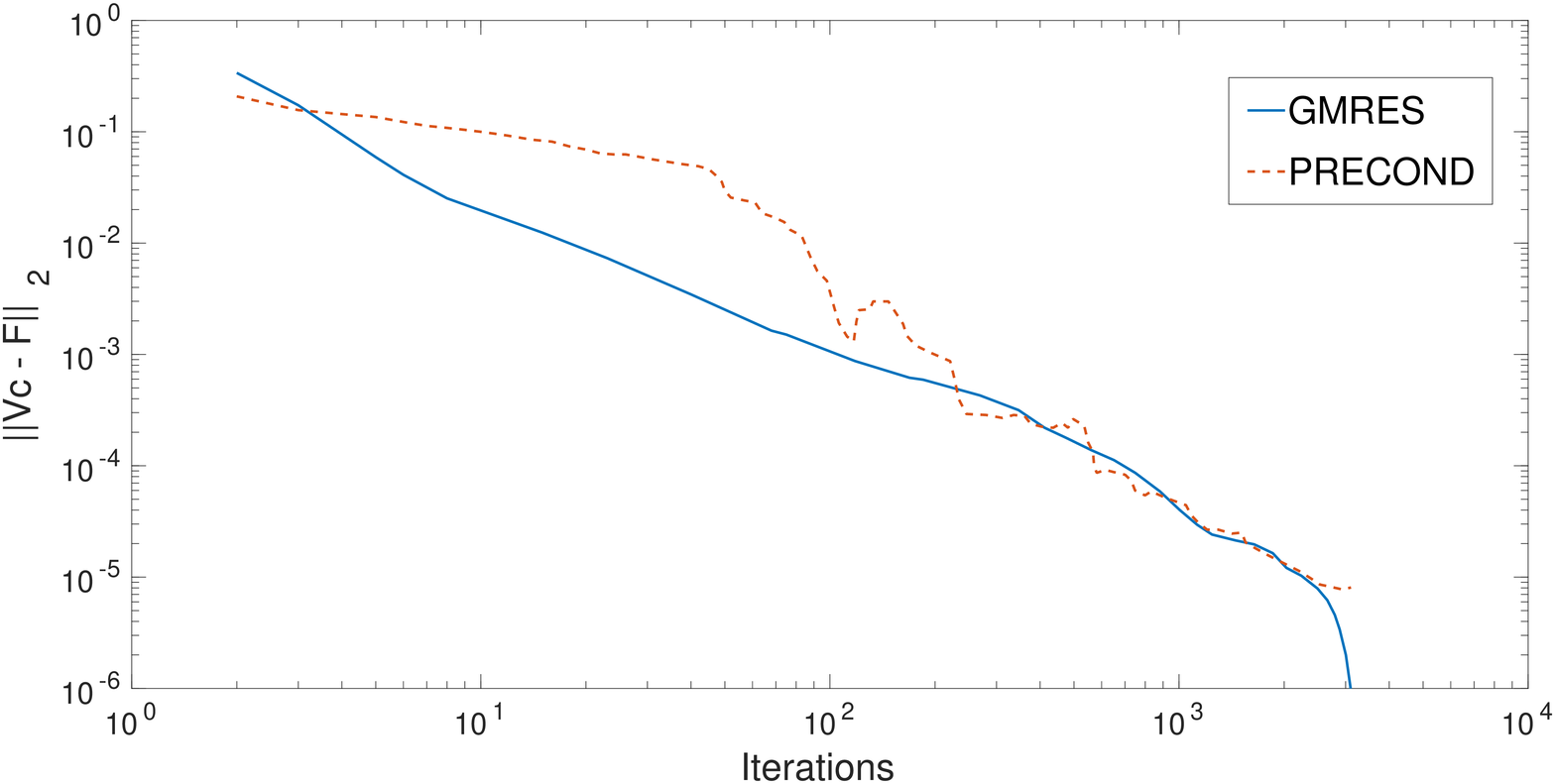}
   \caption{Residual of GMRES for the partition of unity method on the icosahedron with $15$ enrichment functions per triangle.}\label{respu}
\end{figure}

In Figure \ref{precondpu} we plot the number of iterations needed until the standard residual error indicator for GMRES is smaller than $<10^{-6}$ for increasing numbers of enrichments. A comparison between the residuals of GMRES and preconditioned GMRES for $15$ enrichment functions is shown in Figure \ref{respu}.

Even though the preconditioner is not justified here, it still significantly reduces the number of iterations for small numbers of enrichment functions. For larger numbers of enrichment functions, however, the reduction in the final  iteration count becomes less clear. From Table \ref{tab1step}, a simple application of the preconditioner reduces the energy error to $2.1\%$ for $7$, resp.~$1\%$ for $15$ enrichment functions, much below the error of a single preconditioned GMRES iteration.
Note that the partition of unity method leads to poorly-conditioned linear systems. So in addition to using an unjustified preconditioner, for larger degrees of freedom floating point errors become relevant. The condition number here (of the $V^0$ matrix) is up to $3.8\times10^{13}$.

\section{{
Conclusions}}

{
In this work we have presented a preconditioner for the space time systems arising from time domain Galerkin boundary element discretisations. It approximates the algebraic system using extrapolation and may be used as either a preconditioner in a GMRES solver or as a fast, almost exact independent solver of the space-time system for large numbers of degrees of freedom. 
}
{
The numerical experiments show the stability of the preconditioned GMRES  for the standard $h$-TDBEM with linear ansatz and test functions under mesh refinements. The error of the approximate solutions obtained from the preconditioner as an independent solver turns out to be small compared to the discretization error, and it decreases  as the mesh size $h$ tends to $0$, provided the CFL ratio $\frac{h}{\Delta t}$ is fixed. However, the rigorous numerical analysis of the surprisingly good stability properties remains open, while higher-order extrapolation for ansatz and test functions of higher polynomial degree is the content of current work \cite{Gimperleinp}.}

{
For the non-polynomial basis functions of a partition of unity method with a large time step the preconditioner may still reduce the number of necessary GMRES iterations. 
}


\begin{thebibliography}{99}

\bibitem{ajrt} T.~Abboud, P.~Joly, J.~Rodriguez, I.~Terrasse, \emph{Coupling discontinuous Galerkin methods and retarded potentials for transient wave propagation on unbounded domains}, J.~Comp.~Phys.~\textbf{230} (2011), 5877--5907.

\bibitem{aimineu} A.~Aimi, M.~Diligenti, A.~Frangi, C.~Guardasoni, \emph{Neumann exterior wave propagation problems: computational aspects of 3D energetic Galerkin BEM}, Comput.~Mech.~\textbf{51} (2013), 475--493.

\bibitem{aiminumeralg}  A.~Aimi, M.~Diligenti, A.~Frangi, C.~Guardasoni, \emph{A stable 3D energetic Galerkin BEM approach for wave propagation interior problems}, Eng.~Anal.~Bound.~Elem.~\textbf{36} (2012), 1756--1765.

\bibitem{aimi5} A.~Aimi, M.~Diligenti, C.~Guardasoni, I.~Mazzieri, S.~Panizzi, \emph{An energy approach to space-time Galerkin BEM for wave propagation problems}, Internat.~J.~Numer.~Methods Engrg.~\textbf{80} (2009), 1196--1240.


\bibitem{bamberger} A. Bamberger, T. Ha Duong, \emph{Formulation variationnelle espace-temps pour le calcul par potentiel retard de la diffraction d'une onde acoustique}, Math.~Meth.~Appl.~Sciences \textbf{8} (1986), 405--435 and 598--608.

\bibitem{Banz} L.~Banz, H.~Gimperlein, Z.~Nezhi, E.~P.~Stephan, \emph{Time domain BEM for sound radiation of tires}, Comput.~Mech.~\textbf{58} (2016), 45--57.

\bibitem{dfhd} Y.~Ding, A.~Forestier, T.~Ha Duong, \emph{A Galerkin scheme for the time domain integral equation of acoustic scattering from a hard surface}, J.~Acoust.~Soc.~Am.~\textbf{86} (1989), 1566--1572.


\bibitem{survey} H.~Gimperlein, M.~Maischak, E.~P.~Stephan, \emph{Adaptive time domain boundary element methods and engineering applications}, Journal of Integral Equations and Applications \textbf{29} (2017), 75--105.

\bibitem{Gimperlein3} H.~Gimperlein, F.~Meyer, C.~\"{O}zdemir, D.~Stark, E.~P.~Stephan, \emph{Boundary elements with mesh refinements for the wave equation}, Numerische Mathematik \textbf{139} (2018), 867--912.

\bibitem{Nezhi} H.~Gimperlein, Z.~Nezhi, E.~P.~Stephan, \emph{A priori error estimates for a time-dependent boundary element method for the acoustic wave equation in a half-space}, Mathematical Methods in the Applied Sciences \textbf{40} (2017), 448--462.

\bibitem{Gimperlein} H.~Gimperlein, C.~\"{O}zdemir, D.~Stark, E.~P.~Stephan, \emph{A residual a posteriori error estimate for the time domain boundary element method}, preprint.

\bibitem{Gimperleinp} H.~Gimperlein, C.~\"{O}zdemir, D.~Stark, E.~P.~Stephan, \emph{hp-version time domain boundary elements for the wave equation on quasiuniform meshes}, preprint.

\bibitem{Gimperlein2} H.~Gimperlein, C.~\"{O}zdemir, E.~P.~Stephan, \emph{Time domain boundary element methods for the Neumann problem and sound radiation of tires}, J.~Comp.~Math.~\textbf{36} (2018), 70--89.

\bibitem{pupaper} H.~Gimperlein, D.~Stark, \emph{Algorithmic aspects of enriched time domain boundary element methods},
Eng.~Anal.~Bound.~Elem.~(2018), online first, https://doi.org/10.1016/j.enganabound.2018.02.010.


\bibitem{Glafke} M.~Gl\"{a}fke, {\em Adaptive Methods for Time Domain Boundary Integral Equations}. {Ph.D.~thesis, Brunel University, 2012}.

\bibitem{GM1} M.~Gl\"{a}fke, M.~Maischak, {\em Regularity of 2D time domain boundary integral operators}, preprint.

\bibitem{GM2} M.~Gl\"{a}fke, M.~Maischak, {\em Wave-number explicit generalised mapping properties for Helmholtz boundary integral operators, applied to time domain boundary integral operators}, preprint.

\bibitem{Ha-Duong03a} T.~Ha Duong, \emph{On retarded potential boundary integral equations and their discretizations}, in: Topics in computational wave propagation, pp.~301--336, Lect.~Notes Comput.~Sci.~Eng., \textbf{31}, Springer, Berlin, 2003.

\bibitem{hld03} T.~Ha Duong, B.~Ludwig and I.~Terrasse, {\em A Galerkin BEM for transient acoustic scattering by an absorbing obstacle}, Internat.~J.~Numer.~Methods Engrg. \textbf{57} (2003), 1845--1882.

\bibitem{sv} S.~Sauter, A.~Veit, \emph{Adaptive time discretization for retarded potentials}, Numer.~Math.~{\textbf{132} (2016), 569--595}.

\bibitem{sv2} S.~Sauter, A.~Veit, \emph{A Galerkin method for retarded boundary integral equations with smooth and compactly supported temporal basis functions}, Numer.~Math.~\textbf{123} (2013), 145--176. 

\bibitem{sayas} F.-J.~Sayas, \emph{Retarded Potentials and Time Domain Boundary Integral Equations: a road-map}, {Springer Series in Comp.~Math.~50, Springer, 2016}.


\bibitem{terrasse93} I.~Terrasse, \emph{R\'{e}solution math\'{e}matique et num\'{e}rique des \'{e}quations de Maxwell instationnaires par une m\'{e}thode de potentiels retard\'{e}s}, Ph.D.~thesis, \'{E}cole Polytechnique, Palaiseau, 1993.

\bibitem{czech} A.~Veit, M.~Merta, J.~Zapletal, D.~Lukas, \emph{Efficient solution of time-domain boundary integral equations arising in sound-hard scattering}, Internat.~J.~Numer.~Methods Engrg.~\textbf{107} (2016), 430--449.

\bibitem{veit} A.~Veit, \emph{Numerical methods for the time domain boundary integral equations}, Ph.D.~thesis, Universit\"{a}t Z\"{u}rich, 2011.

\bibitem{mic} A.~E.~Yilmaz, J.-M.~Jin, E.~Michielssen, \emph{Time domain adaptive integral method for surface integral equations}, IEEE Trans.~Antennas Propagation \textbf{52} (2004) 2692--2708.


\end{thebibliography}
\end{document}